\documentclass[11pt]{article}
\usepackage[affil-it]{authblk}
\usepackage{mathrsfs}

\newtheorem{rem}{Remark}

\usepackage{latexsym}
\usepackage{amsmath,mathtools}
\usepackage{amssymb}
\usepackage{graphicx}
\usepackage{textcomp}
\usepackage{multirow}
\usepackage{enumerate}
\usepackage{cancel}
\usepackage{caption}  
\usepackage[usenames,dvipsnames]{color}
\usepackage[left, pagewise]{lineno}
\usepackage[ruled,vlined]{algorithm2e}
\usepackage{epstopdf}
\usepackage{hyperref}
\newcommand{\wh}{\widehat}

\newcommand{\bF}{{\bf F}}
\newcommand{\bI}{{\bf I}}
\newcommand{\bP}{{\bf P}}

\newcommand{\bU}{{\bf U}}
\newcommand{\bX}{{\bf X}}

\newcommand{\bPhi}{\boldsymbol{\Phi}}
\newcommand{\bPsi}{\boldsymbol{\Psi}}

\title{A FOM/ROM Hybrid Approach for Accelerating Numerical Simulations 
\thanks{The authors acknowledge the support from the Institute for Computational and Experimental Research in Mathematics (ICERM) at Brown University for participating the ``Model and dimension reduction in uncertain and dynamic systems" semester program in Spring 2020, which were supported by the National Science Foundation under Grant No. DMS-1439786 and the Simons Foundation Grant No. 50736.
Z. Wang was partially supported by the National Science Foundation through Grants No. DMS-1913073,2012469 and the U.S. Department of Energy Grant No. DE-SC0020270. 
}
}

\author{Lihong Feng
	\thanks{E-mail: \texttt{feng@mpi-magdeburg.mpg.de}}} 
	\affil{Max Planck Institute for Dynamics of Complex Technical Systems, 1 Sandtorstr, 39106 Magdeburg, Germany}
	
\author{Guosheng Fu
	\thanks{E-mail: \texttt{gfu@nd.edu}}}
	\affil{ Department of Applied and Computational Mathematics and Statistics, University of Notre Dame, Notre Dame, IN 46556}
	
\author{Zhu Wang
	\thanks{E-mail: \texttt{wangzhu@math.sc.edu}; Corresponding author}}
	\affil{ Department of Mathematics, University of South Carolina, Columbia, SC 29208}

\date{}
\begin{document}
\maketitle
\begin{abstract} 
The basis generation in reduced order modeling usually requires multiple high-fidelity large-scale simulations that could take a huge computational cost. In order to accelerate these numerical simulations, we introduce a FOM/ROM hybrid approach in this paper. It is developed based on an a posteriori error estimation for the output approximation of the dynamical system. By controlling the estimated error, the method dynamically switches between the full-order model and the reduced-oder model generated on the fly. Therefore, it reduces the computational cost of a high-fidelity simulation while achieving a prescribed accuracy level. Numerical tests on the non-parametric and parametric PDEs illustrate the efficacy of the proposed approach. 
\end{abstract}

\section{Introduction}
\label{sec:in}
Reduced order modeling has been developed to provide an efficient computational model for dynamical systems appeared in many scientific research and engineering application problems. 
Given a parametric, time-dependent problem, the projection-based model reduction method first generates reduced basis at an offline stage from snapshots obtained at selected parameters and time instances, and builds a low-dimensional model in which the matrices and vectors of the system can be fully assembled; then at the online stage, this reduced-order model (ROM) can be simulated at a very low cost, thus make fast or even real-time numerical simulations feasible for many-query applications. 
Model reduction techniques have achieved many successes, but also faces several issues.  

On the one hand, to build ROMs for a complex system, one usually has to run large-scale, high fidelity simulations for many times to find the reduced basis. Although the basis generation is part of the offline process, its high computational effort makes reduced order modeling not attractive in practice. Thus, it is desired to speed up this offline process. 
On the other hand, it is difficult to build effective ROMs for certain problems. 
For instance, if problems with a large Kromogrov n-width are considered such as transport or convection-dominated phenomena, it is generally hard to find a low-dimensional basis for approximating the solution manifold. 
Meanwhile, problems with varying initial conditions or forcing can bring big challenges to the ROM as well because it is impossible to parametrize the general initial condition or the forcing, and snapshots are not likely to cover all the possible system responses. 
This also makes the prediction ability of ROMs pretty weak: if the snapshots are collected over a given time interval, once the ROM evolves beyond the interval, there is no accuracy guarantee of the reduced-order simulations. 

Several approaches have been pursued to tackle these issues. 
To save the computational cost for snapshot generation, a hybrid snapshot simulation methodology was proposed to accelerate the high-quality data generation for the proper orthogonal decomposition (POD) ROM development in \cite{bai2021reduced,bai2020deim}.  In this approach, the simulation alternates between full-order model (FOM) and local ROM based on criteria related to the singular value distribution and the magnitude of POD basis coefficients. 
Similarly, the equation-free/Galerkin-free approaches developed in \cite{esfahanian2009equation,Sirisup2005A} also chose between the fine-scale and coarse-scale models during the simulations. The fine-scale run provides data by performing full-order simulations with a small time step size that finds a handful of leading POD basis. The POD modes are used to parametrize the low-dimensional attracting slow manifold; while the coarse-scale run evolves the coefficients of leading POD basis functions by using a larger step size. The resulting reduced-order approximation can be constructed accordingly, which then initializes the fine-scale simulations. 
However, the choice of the frequency for switching between the two scale models is empirical. 
Note that the idea combining the ROM and FOM has been used in parallel-in-time simulations, for example, in  \cite{farhat2006time,chen2014use,maday2020adaptive}, where the ROM on the fly is regarded as a coarse propagator. 
For problems with large Kromogrov n-width, the idea of adaptivity has been introduced for adjusting reduced basis dynamically. The adaptive h-refined POD was developed in \cite{carlberg2015adaptive,etter2020online}. In this approach, an adjoint-based error estimation is used to mark basis vectors with major error contributions. Such basis vectors are then refined till the estimated error becomes under a prescribed tolerance. The refined basis will be compressed online if the refined-basis dimension exceeds a specified threshold or a prescribed number of time steps has elapsed. An online adaptive bases and adaptive sampling was developed for transport-dominated problems in \cite{peherstorfer2020model}, in which both the state and the nonlinear term are projected onto the POD basis subspace, and the existing POD basis is adjusted at every step. The basis update is of low rank that is sought from an optimization problem for minimizing the POD/DEIM projection error of the nonlinear term. A related approach has been proposed in \cite{hesthaven2020rank} for Hamiltonian systems. 
Another research direction is to consider the data-driven ROMs, which first finds nonlinear solution manifold by taking advantages of the expression power of deep neural network, and then builds projection-based or machine learning-based ROMs in the nonlinear manifold, for instance, in \cite{lee2020model,lee2019deep,kim2020efficient,fresca2021pod,san2018neural,maulik2021reduced}.

In this paper, we focus on the projection-based ROMs and use the POD method for constructing the reduced basis. In order to accelerate the snapshot generation and deal with situations for which the ROM may be not effective, we propose a FOM/ROM hybrid approach. The main idea is to dynamically update ROM on the fly and switch simulation between the FOM and the ROM, which is close to the idea in \cite{esfahanian2009equation,Sirisup2005A,bai2021reduced}, but we introduce a more rigorous criterion based on an {\em a posteriori} error indicator. By comparing the estimated error with a user-defined tolerance, the hybrid approach automatically alternates between FOM and ROM while controlling the output approximation error.  

The rest of the paper is organized as follows. In Section \ref{sec:err}, we discuss the {\em a posteriori} error indicator; in Section \ref{sec:alg}, we introduce the hybrid approach, which is numerically investigated in Section \ref{sec:num}. A few concluding remarks are drawn in the last section. 

\section{Error estimation of output approximation}
\label{sec:err}
  
Consider the following discrete input/output dynamical system: given the state $x^k$ and input $u^k$ at time $t_k$, to find the output $y^{k+1}$ that satisfies 
\begin{subequations}
\begin{align}%
E^k x^{k+1} &= A^k x^{k} + f(x^k) + B^k u^k, \label{eq:fom01} \\
y^{k+1} &= Cx^{k+1}. \label{eq:fom02}
\end{align}
\label{eq:fom0}
\end{subequations}
The system in \eqref{eq:fom0} is regarded as the full-order model (FOM). The first equation of the state may come from a time-dependent PDE after discretizing the differential operators and state variable in space and considering an explicit/semi-implicit time-stepping algorithm. In case an implicit scheme is used for integrating a nonlinear PDE, this equation could appear in the linearized system during nonlinear iterations. 
The second equation evaluates the output of the system, which approximates linear functionals of the PDE solution -- quantities of interest of the system. 
Naming \eqref{eq:fom01} the {\em primal} equation and, 
to measure and further control the output approximation error, we follow \cite{chellappa2020adaptive} and consider the following {\em dual} problem: to find $x_{du}^{k+1}$ at time $t_{k+1}$ satisfying
\begin{equation}
(E^{k})^\intercal x_{du}^{k+1} = -C^\intercal. 
\label{eq:dual}
\end{equation}
Based on the original system and choices of time step sizes, the coefficient matrices $E^k, A^k$ and $B^k$ could vary with time. 
Here, for simplicity of presentation, we assume that they stay unchanged during the simulation in the sequel of discussion. Therefore, we instead consider the following primal equation: 
\begin{equation}
E x^{k+1} = A x^{k} + f(x^k) + B u^k,
\label{eq:fom1}
\end{equation}
together with the output equation \eqref{eq:fom02} and the dual equation: 
\begin{equation}
E^\intercal x_{du} = -C^\intercal.
\label{eq:dual1}
\end{equation} 
But we emphasize that our approach could be directly applied to the general case represented by \eqref{eq:fom0} and \eqref{eq:dual}. 

To reduce the computational cost of the FOM simulations, model reduction techniques such as the POD method can be used to generate a ROM. 
The POD method extracts a low-dimensional set of orthonormal basis vectors from the snapshot data and uses it to build a low-dimensional model based on projections.  
Indeed, given the POD basis $\Phi_k$ and the state $x^k$ at time $t_k$, we find a reduced-order state approximation $\wh{x}^{k+1}= \Phi_k x_r^{k+1}$ to approximate $x^{k+1}$ and further determine a reduced-order output approximation $y_r^{k+1}$. The associated Galerkin projection-based ROM has the following form: 
\begin{subequations}
\begin{align}
E_r  x_r^{k+1} &= A_r x_r^{k} + f_r(\Phi_k x_r^k) + B_r u^k, \label{eq:rom_p}\\
y_r^{k+1} &= C_r x_r^{k+1}, \label{eq:rom_d}
\end{align}
\label{eq:rom1}
\end{subequations}
where 
$x_r^k = \Phi_k^\intercal x^k$, 
$E_r=\Phi_k^\intercal E \Phi_k$, 
$A_r =\Phi_k^\intercal A \Phi_k$,   
$B_r = \Phi_k^\intercal B$, 
$C_r = C \Phi_k$, 
and 
$f_r(\Phi_k x_r^k) = \Phi_k^\intercal \wh{f}(\Phi_k x_r^k)$ with 
$\wh{f}(\cdot) = f(\cdot)$ if POD is used or 
$\wh{f}(\cdot) = \mathbb{P}f(\cdot)$ if an interpolation method is applied. For instance, when DEIM is considered, $\mathbb{P} = \bPhi_f(\bP^\intercal \bPhi_f)^{-1}\bP$ with $\bPhi_f$ the DEIM basis and $\bP$ the matrix for extracting rows corresponding to the interpolation points. 

The dual equation \eqref{eq:dual1} needs to be considered for estimating the error, which has the same dimension as the FOM. To avoid solving the large-scale discrete system, one can seek an approximate solution $\widehat{x}_{du}$, e.g., via a ROM of the dual equation. Let $\widehat{x}_{du} = \bPsi_k x_r^{du}$ with $\bPsi_k$ the POD basis of $x_{du}$, the ROM of dual equation reads:  
\begin{equation}
E_r^\intercal x_r^{du} = -\bPsi_k^\intercal C^\intercal. 
\label{eq:dual1b}
\end{equation}

To derive a bound on the output error, we test \eqref{eq:dual1} by $x^{k+1}-\wh{x}^{k+1}$ and have  
\begin{equation}
(x^{k+1}-\wh{x}^{k+1})^\intercal E^\intercal x_{du} = -(x^{k+1}-\wh{x}^{k+1})^\intercal C^\intercal. 
\label{eq:dual2}
\end{equation}
Taking transpose on both sides yields 
\begin{equation}
C (x^{k+1}-\wh{x}^{k+1}) = - x_{du}^\intercal E (x^{k+1}-\wh{x}^{k+1}), 
\label{eq:dual3}
\end{equation}
in which the term on the left represents the output error of the reduced approximation and the term on right is related to the solution of the dual problem and the residual of the primal problem. 
Define 
\begin{equation}
r_{pr}^{k+1} \coloneqq A x^{k} + f(x^k) + Bu^k - E \wh{x}^{k+1}
\overset{\eqref{eq:fom1}}{=} E (x^{k+1}-\wh{x}^{k+1}), 
\label{eq:dual4}
\end{equation}
and define the residual of dual equation to be 
\begin{equation}
r_{du} \coloneqq -C^\intercal - E^\intercal \wh{x}_{du}
\overset{\eqref{eq:dual1}}{=} E^\intercal (x_{du}-\wh{x}_{du}).   
\label{eq:dual5}
\end{equation}
Therefore, combining \eqref{eq:dual3}, \eqref{eq:dual4} and \eqref{eq:dual5}, we have 
\begin{equation}
\begin{aligned}
y^{k+1}-y_r^{k+1} &= C (x^{k+1}-\wh{x}^{k+1}) = - x_{du}^\intercal r_{pr}^{k+1} \\
		            &= [- (x_{du}^\intercal-\wh{x}_{du}^\intercal) - \wh{x}_{du}^\intercal] r_{pr}^{k+1} \\
		            &= -( r_{du}^\intercal E^{-1}+ \wh{x}_{du}^\intercal ) r_{pr}^{k+1} .       
\end{aligned}
\label{eq:outputerr}
\end{equation}
It leads to an upper bound for the output approximation: 
\begin{equation}
|y^{k+1}-y_r^{k+1}| 
		\leq \left( \| r_{du} \| \|E^{-1}\|  + \|\wh{x}_{du}\| \right) \|r_{pr}^{k+1}\|.    
\label{eq:outputerr2}
\end{equation}
However, evaluating $r_{pr}^{k+1}$ in \eqref{eq:dual4} needs the FOM solution $x^k$. To reduce the computational cost, we replace $r_{pr}^{k+1}$ with the residual of the reduced-order approximation, i.e., 
\begin{equation}
\wh{r}_{pr}^{k+1} \coloneqq A \wh{x}^{k} + f(\wh{x}^k) + Bu^k - E \wh{x}^{k+1}.  
\label{eq:dualres}
\end{equation}
Then the output error satisfies 
\begin{equation}
|y^{k+1}-y_r^{k+1}| 
		\leq \left( \| r_{du} \| \|E^{-1}\|  + \|\wh{x}_{du}\| \right) \rho^{k+1} \|\wh{r}_{pr}^{k+1}\|,    
\label{eq:outputerr3}
\end{equation}
where 
$$\rho^{k+1} = \frac{\|{r}_{pr}^{k+1}\|}{\|\wh{r}_{pr}^{k+1}\|}.$$
It is proved in \cite{zhang2015efficient} that, under mild conditions, e.g., $f$ is bi-Lipschitz continuous, the ratio $\rho^{k+1}$ is bounded by positive constants from above and below. 
In particular, 
$$
\|r_{pr}^{k+1}-\wh{r}_{pr}^{k+1}\| 
= \| A (x^k - \wh{x}^{k}) + f(x^k) - f(\wh{x}^k)\|
\leq (\|A\|+L_f) \|x^k - \wh{x}^{k}\|.  
$$
When the reduced approximation is accurate enough, i.e., $\|\wh{x}^k-x^k\|$ is sufficient small, $\|r_{pr}^{k+1}-\wh{r}_{pr}^{k+1}\|$ becomes close to zero, and the ratio $\rho^{k+1}$ approaches 1. 
Indeed, this is true in our computational setting presented in the next section, as we aim at finding accurate numerical solutions by combining the FOM and the ROM during the simulations.  
Therefore, we define 
\begin{equation}
\Delta^{k+1} = \left( \| r_{du} \| \|E^{-1}\|  + \|\wh{x}_{du} \| \right) \|\wh{r}_{pr}^{k+1}\|,
\label{eq:est}
\end{equation}
which gives an error estimation of $|y^{k+1}-y_r^{k+1}|$ for the system \eqref{eq:fom1}. 
We emphasize that $\Delta^{k+1}$ only serves as an error indicator at the step $k$ because the numerical error accumulates over time that may deflect ${x}^k$, the initial condition of the very step, off the trajectory of a high-fidelity simulation. 

\section{A hybrid approach}
\label{sec:alg}

Based on the error indicator presented in Section \ref{sec:err}, we introduce a FOM/ROM hybrid approach in this section. 
The goal is to accelerate numerical simulations while keeping an accurate output approximations. 
The basic idea of this approach is to generate a ROM on the fly based on snapshots available in a fixed-width time window, and apply the ROM whenever possible. The associated output accuracy is measured by the error indicator $\Delta^k$. If it does not reach a user-defined tolerance, one switches the model to the FOM, updates snapshots, and generates the ROM that is to be used in the next step's simulation. 
Therefore, this approach dynamically switches the simulation model between the ROM and the FOM. 
For the system, \eqref{eq:fom1} and \eqref{eq:dual1}, the approach is detailed in Algorithm \ref{alg:hyb}.

\begin{algorithm}[!ht]\footnotesize
\SetAlgoLined
\KwIn{ Initial condition $x_0$, error tolerance $\mathsf{tol}$, time window width $w$, total number of time steps $n_t$.}
\KwOut{ System output.}
compute an approximate solution $\widehat{x}_{du}$ to the dual equation and evaluate $C_{du}= \| r_{du} \| \|E^{-1}\|  + \|\wh{x}_{du}\|$\;
run the FOM from $t^1$ to $t^w$ and compute the output; generate snapshots $S=[x^1, x^2, \ldots, x^w]$, $S_E = ES$, $S_A = AS$ and initialize a ROM\;
\For{$k=w:n_t-1$}{
	run the ROM on $[t^k, t^{k+1}]$, obtain $\wh{x}^{k+1}$ and evaluate the estimated error $\Delta^{k+1}= C_{du}\|\widehat{r}_{pr}^{k+1}\|$\;
	\eIf{$\Delta^{k+1}\leq \mathsf{tol}$}{
	$x^{k+1} = \wh{x}^{k+1}$, $\mathsf{flag}= 1$\; }{
	{
	solve the FOM for $x^{k+1}$\; 
	set $\Delta^{k+1}= \mathsf{eps}$ (in MATLAB notation) and $\mathsf{flag}= 2$.}
	}
	evaluate output $y^k$, update time window data: $S$, $S_E$ and $S_A$ \; 
 
	\If{$\mathsf{flag}=2$}{generate a ROM based on $S$, $S_E$ and $S_A$.}
}
\caption{The hybrid approach for solving \eqref{eq:fom1}-\eqref{eq:fom02} \label{alg:hyb}}
\end{algorithm}


Note that here we still assume time invariant coefficient matrices in system \eqref{eq:fom1}. 
Since the dual problem \eqref{eq:fom02} only needs to be solved once, we use a full-order approximation and compute $C_{du} = \| r_{du} \| \|E^{-1}\|  + \|\wh{x}_{du}\|$ at the beginning of the simulation. 
The hybrid approach starts with the FOM simulation over a fixed-width time window, 
$[t^1, \ldots, t^w]$. Based on this time window data $S=[x^1, x^2, \ldots, x^w]$, a ROM is initialized; 
At each subsequent time step $k$, we first use the available ROM to advance the simulation over a step, and evaluate the output error indicator. By comparing it with a prescribed tolerance ($\mathsf{tol}$), we will either accept the reduced approximation $\wh{x}^{k+1}$ as $x^{k+1}$ if $\Delta^{k+1}\leq \mathsf{tol}$, or reject it if $\Delta^{k+1}> \mathsf{tol}$ and run a full-order simulation to get $x^{k+1}$ with a higher fidelity.  
In either case, $x^{k+1}$ is used to update the time window data $S$.   
Whenever a FOM solution is added, we shall update the ROM.


The computational saving stems from the replacement of FOM with ROM at all the possible time steps. 
To determine when to use FOM, the error indicator $\Delta^{k+1}$ is evaluated where the major computation is spent on calculating $\|\wh{r}_{pr}^{k+1}\|$. 
For updating the ROMs, the main computation is spent on updating the POD basis $\Phi_{k+1}$ from $S$ and assembling coefficient matrices and vectors, $A_r$, $E_r$ and $f_r$, in the reduced system.  
Note that 
$$
\wh{r}_{pr}^{k+1} = A\Phi_k x_r^{k} + f(\Phi_k x_r^{k}) + Bu^k - E \Phi_k x_r^{k+1}, 
$$
and 
$A_r =\Phi_k^\intercal A \Phi_k$, 
$E_r=\Phi_k^\intercal E \Phi_k$,  
calculating $A\Phi_k$ and $E\Phi_k$ in a cheaper way would improve the efficiency. 
On the other hand, the POD basis $\Phi_k$ is generated based on the available time window snapshots $S$, which is a tall matrix in general. 
For such a case, the method of snapshots (MOS) can be used to find the POD basis (i.e., left singular value vectors of $S$) quicker than the singular value decomposition (SVD): 
$$
S^\intercal S V = V \Lambda \text{ and } \Phi = S V\Lambda^{-1/2}.
$$
The method can be further improved in efficiency as discussed in \cite{wang2016approximate} by taking advantage of parallel computing. 
Thus, at the end of the step $k$ in the algorithm, besides updating the time window snapshots $S=[S(:,\,2:w),x^{k+1}]$, we also update $S_E= [S_E(:,\, 2:w), Ex^{k+1}]$ and $S_A = [S_A(:,\,2:w), Ax^{k+1}]$. 
It leads to a cheaper evaluation of 
$$E\Phi_{k+1} = S_E V \Lambda^{-1/2} \text{ and } A\Phi_{k+1} = S_A V \Lambda^{-1/2}$$
 at the next step.  

The computational complexity at each step is dominated by the time window data updates if $\mathsf{flag}=1$. 
There, the main complexity, $\mathcal{O}(n\max(r,a))$ with $r$ the dimension of ROM and $a$ the max nonzero numbers in each row of sparse matrices $E$ and $A$, is used to update time window data and evaluate error indicator. 
It generally dominates the cost for solving the reduced system. 
When $\mathsf{flag}=2$, the computational complexity at each step is dominated by the FOM solution and the ROM update. 
There, the FOM solution takes $\mathcal{O}(n^p)$ with $p$ dependent on the choice of numerical solver, suppose it dominates the nonlinear function evaluation cost $\mathcal{O}(\alpha(n))$, and the MOS takes $\mathcal{O}(nw^2+rnw+w^3)$, in which $\mathcal{O}(nw^2)$ flops are needed for evaluating $S^\intercal S$, $\mathcal{O}(w^3)$ flops for the eigen-decomposition of $S^\intercal S$ and $\mathcal{O}(nwr)$ flops for the calculation of the POD basis. 
Note that, when the DEIM is used, additional computation is required for updating the DEIM basis and interpolation points. 
Since $r\leq w$ and $w\ll n$ in general, the dominating complexity becomes $\mathcal{O}(nw^2)$. 
Assume the percentage of the FOM runs in the hybrid model simulation is $c\times 100\%$, for $c\in [0, 1]$, then the ratio of the complexity for the FOM to that of the hybrid approach is $\frac{n^p}{c(n^p+nw^2)+(1-c)(n\max(r,a))}$. Obviously, this factor is bounded by $\frac{1}{c}$ from above and by $1$ from below when $w^2\leq (1-c)n^{(p-1)}$. 

The hybrid approach accelerates the FOM simulations, which has several applications:
Firstly, it can be used to speed up the offline snapshot generation for model order reduction (MOR) techniques as it replaces the FOM simulations with the ROM simulations at certain steps; 
Secondly, it can overcome the issue on the accuracy loss when using a ROM in a computational setting different from which the ROM was built. 
Next, we focus on the first application to illustrate the usefulness of the approach. 
\section{Accelerating the offline stage of MOR}
\label{sec:acc}
The application of ROM usually involves certain high-fidelity simulations to find reduced-order basis vectors offline and/or update them at the online stages. When the system has a large dimension, generation of snapshots using the FOM could be time-consuming. 
In this section we show that the hybrid approach can be naturally employed at the offline stage of MOR.

In particular, when the reduced basis method (RBM) is applied to parametric time-evolution problems, the basis to approximate the solution manifold can be found by using the POD-Greedy algorithm \cite{haasdonk2008reduced}. Guided by an a posteriori error estimator $\Delta(\cdot)$, the algorithm implements a sequence of FOM simulations at selected parameter samples to obtain the corresponding trajectory of the solution, i.e. the snapshots. At this point, the hybrid simulation in Section \ref{sec:alg} can be applied to accelerating the FOM simulation with accuracy loss being controlled by the error indicator. 
In Algorithm \ref{alg:adp0}, we present a modified, hybrid model-based POD-greedy algorithm. 
It is worth mentioning that, because the FOM runs occur at a portion of time steps, e.g., $\{t_1, \ldots, t_K\}$ in the hybrid model simulations at $\mu^*$ and the ROMs on the fly are essentially generated from them, one only needs to collect the snapshots $\bX_{\mu^*}$ at those steps, that are enough to represent the system information. 
\begin{algorithm}[!htb]\footnotesize
\SetAlgoLined
\KwIn{Parameter training set $\Theta$, tolerance $\mathsf{tol}$.}
\KwOut{basis matrix $\bPhi$.}
Initialize $\bPhi=[\,\,]$, randomly choose $\mu^*\in \Theta$\;
\While{$\Delta(\mu^*)> \mathsf{tol}$}{
	simulate the {\bf hybrid} model at $\mu^*$ and obtain snapshots $\bX_{\mu^*}=[x(t_1, \mu^*), \ldots, x(t_K, \mu^*)]$\;
	find left singular vectors $\bU$ of $\bX_{\mu^*}= (\bI-\bPhi\bPhi^\intercal) \bX_{\mu^*}$ and define $\bPhi= \mathsf{orth}\{\bPhi, \bU(:,1)\}$\;
	find $\mu^*= \arg\max_{\mu\in\Theta} \Delta (\mu)$.
	}
\caption{The hybrid model-based POD-Greedy algorithm\label{alg:adp0}}
\end{algorithm}

Based on the same idea, we can further apply the proposed hybrid method to the adaptive basis construction approaches introduced in \cite{chellappa2020adaptive}. 
Especially, the adaptive POD-Greedy-DEIM algorithm in \cite{chellappa2020adaptive} builds a nonlinear ROM offline, which chooses the POD basis vectors, the DEIM basis vectors and interpolation points in an adaptive manner controlled by an a posteriori error estimator $\overline{\Delta}(\cdot)$. 
Since high-fidelity simulations need to be performed at selected parameters for updating basis, we replace the associated FOM simulations with the FOM/ROM hybrid ones. 
The modified, hybrid model-based algorithm is presented in Algorithm \ref{alg:adp}. 

\begin{algorithm}[!htb]\footnotesize
\SetAlgoLined
\KwIn{Parameter training set $\Theta$, tolerance $\mathsf{tol}$.}
\KwOut{POD basis $\bPhi$, DEIM basis $\bPhi_f$ and interpolation points $\bP$.}
solve the nonparametric dual system and obtain $C_{du}$\;
Initialize $\bPhi=[\,\,]$, $I_{\mu} = [\,\,]$, $\ell_{POD}=1$, $\ell_{DEIM}=1$, randomly choose $\mu^*\in \Theta$\;
\While{$\overline{\Delta}(\mu^*)\notin \mathsf{zoa}=[\mathsf{tol}/10, \mathsf{tol}]$}{
	\eIf{$\ell_{POD}<0$}{
	remove the last $\ell_{RB}$ columns from $\bPhi$\;
	}{
	\eIf{$\mu^*\notin I_{\mu}$}{
	solve {\bf hybrid} model at $\mu^*$ for snapshots $\bX_{\mu^*}$ and nonlinear snapshots $\bF_{\mu^*}$ and store them\;
	}{
	load snapshots $\bX_{\mu^*}$\;
	}
	find left singular vectors $\bU$ of $\bX_{\mu^*}= (\bI-\bPhi\bPhi^\intercal) \bX_{\mu^*}$ and define $\bPhi= \mathsf{orth}\{\bPhi, \bU(:,1:\ell_{POD})\}$\;
	}
	update $\ell_{DEIM}$ DEIM basis $\bPhi_f$ and interpolation points $\bP$ based on all the available nonlinear snapshots\;
	update the POD-DEIM model and evaluate $\overline{\Delta} (\mu) = \overline{\Delta}_{POD} (\mu) + \overline{\Delta}_{DEIM} (\mu) $ for all $\mu\in \Theta$\;
	find $\mu^*= \arg\max_{\mu\in\Theta} \overline{\Delta} (\mu)$ and update $I_{\mu} = [ I_{\mu}, \mu^*]$\;
	$\ell_{POD} = \ell_{POD}+\lfloor \log\left(\frac{\overline{\Delta}{POD}(\mu^*)}{tol}\right) \rfloor \pm 1$, and 
	$\ell_{DEIM} = \ell_{DEIM}+\lfloor \log\left(\frac{\overline{\Delta}{DEIM}(\mu^*)}{tol}\right) \rfloor \pm 1$.
}
\caption{The hybrid model-based adaptive POD-Greedy-DEIM algorithm\label{alg:adp}}
\end{algorithm}
\section{Numerical results}
\label{sec:num}
To investigate the performance of the proposed hybrid method, 
we consider three test cases in this section: the first one is the non-parametric Burger's equation, the second one is the parametric Burgers' equation, and the last one is a nonlinear circuit problem with a step input signal. 
For the first and third cases, we aim at achieving accurate numerical simulations; 
for the second one, we combine the hybrid method with adaptive offline basis selection for reduced order modeling.

To check the behavior of the error indicator, we compare the estimated error $\Delta^{k+1}$ with the ``true" error that is referred to the value of $|y^{k+1}-y_r^{k+1}|$ for the given data $x^k$ at the time step $k$. 
For investigating the numerical performance of the hybrid approach, we quantify the accuracy of its numerical simulation by $\mathcal{E}_o$, the mean error of the approximate output over the time interval; and measure the efficiency by $P_f$, the percentage of the number of FOM runs in the total number of time steps for the hybrid simulation, and by the wall-clock time for integrating the hybrid model, $t_{h}$, compared to that of the FOM, $t_{f}$. All simulations are implemented in Matlab and the wall-clock time are estimated by the timing functions $\mathsf{tic/toc}$.  
 
\paragraph{Test case 1.} We first consider the 1D Burgers' equation:  
\begin{equation*}
\partial_t q(x, t) -\nu \partial_x^2 q(x, t) + q(x, t) \partial_x q(x, t) = 0
\end{equation*}
for $x\in [0, 1]$ and $t\in [0, 1]$, where the viscosity coefficient $\nu$ is a constant. 
The problem is associated with the zero Dirichlet boundary condition imposed at both endpoints of the domain and the initial condition 
\begin{equation*}
q(x, 0) = \left\{
\begin{array}{ll}
1, &\,\,0.1\leq x\leq 0.2,\\
0, &\,\,\text{otherwise},
\end{array}
\right.
\end{equation*}
together with an output of interest $y(t) = \int_{0}^1 q(x,t)\, dx$. 

When a full order simulation is performed, the semi-implicit Euler method is taken for time stepping. The whole domain is equally partitioned by $n=2^p-1$ interior grid points and the time interval is divided into $N_t=2^{p+1}$ uniform subintervals. 
The finite difference discretization using central scheme leads to the following system:  
\begin{equation*}
E Q^{k+1} = Q^{k}  + f(Q^k), 
\end{equation*}
where $E = I - \Delta t \nu L$ with $L$ the discrete Laplacian and $f(Q^k) = -\Delta t Q^k\circ (D_x Q^k)$ with $D_x$ the discrete first-order differential operator and $\circ$ represents the Hadamard product.   
The output quantity represents the average value of the state variable over the domain 
\begin{equation*}
Y^{k+1} = C Q^{k+1}, 
\end{equation*}
where $C= \frac{1}{n} (\mathbf{1}_n)^\intercal$ with $\mathbf{1}_n$ denotes the all-ones column vector of dimension $n$. 

\begin{figure}[!htb]
 \centering
        \includegraphics[width=1\linewidth,height=2in]{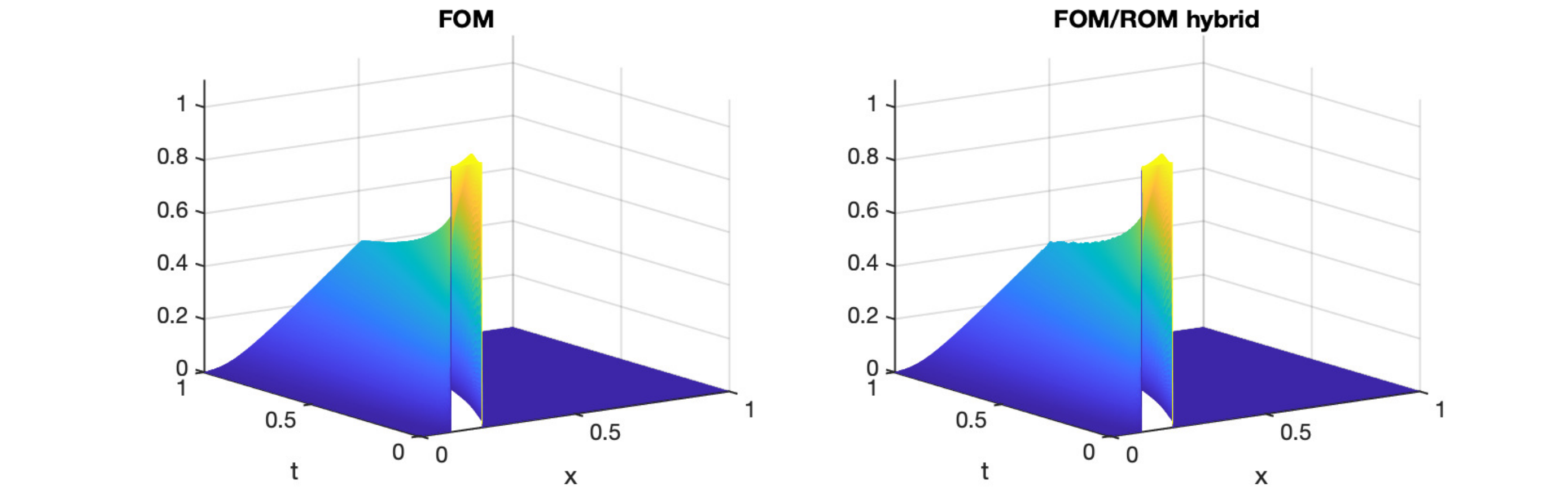}
\caption{Time evolution plots of Burgers' equation at $\nu=10^{-3}$: (left) FOM, (right) the hybrid model when $w=20$ and $tol=10^{-4}$. \label{fig:evo1}}
\end{figure}

\begin{figure}[!htb]
 \centering
        \includegraphics[width=1\linewidth]{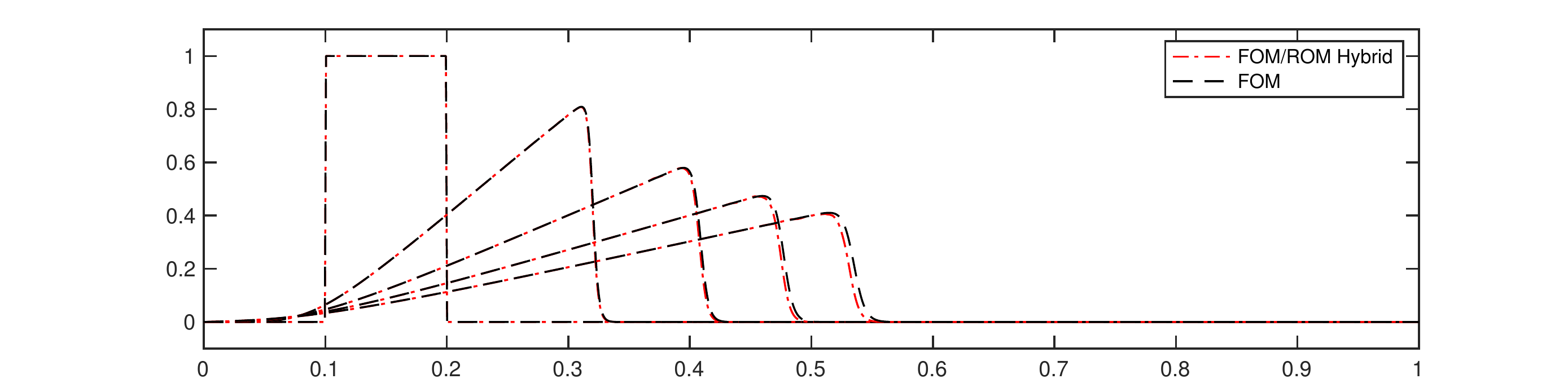}
\caption{Time snapshots $t=0, 0.25, 0.5, 0.75, 1$ of Burgers' equation at $\nu=10^{-3}$, $w=20$. \label{fig:snap1}}
\end{figure}

\begin{figure}[!htb]
 \centering
        \includegraphics[width=1\linewidth]{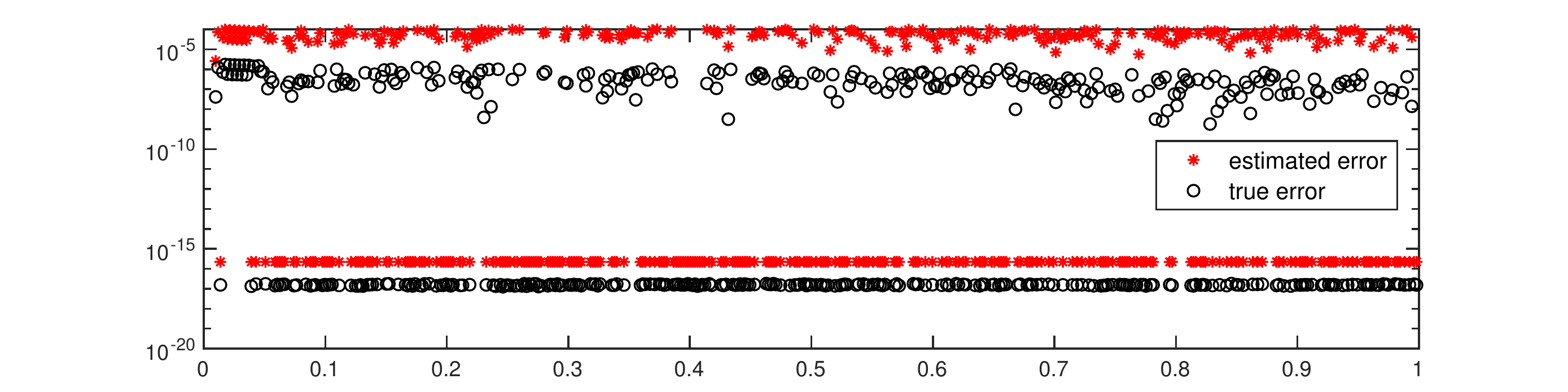}
\caption{Estimated and true errors of Burgers' equation at $\nu=10^{-3}$, $w=20$.\label{fig:err1}}
\end{figure}

Firstly, we set $\nu=10^{-3}$ and $p=10$ in the discretization. 
In the hybrid approach, we use POD ROM, set the time window width $w=20$ and the tolerance of output error at each step to be $10^{-4}$. 
As the dual problem \eqref{eq:dual1} only needs to be solved once, we compute it by a direct solver. 
The time evolutions of the FOM and the hybrid model results are plotted in Figure \ref{fig:evo1}, respectively. 
The snapshots at $t=0, 0.25, 0.5, 0.75, 1$ obtained by both approaches are plotted in Figure \ref{fig:snap1}. 
It is seen that, compared with the FOM, the hybrid approach obtains an accurate state approximation. 
Time evolutions of estimated errors $\Delta^{k+1}$ and ``true" errors, plotted once every 4 steps, are shown in Figure \ref{fig:err1}. 
Note that for some steps, the estimated error is of machine precision that is because full order simulations are performed at those time steps. 
It shows that $\Delta^{k+1}$ provides a reliable indicator of the output approximation error. 


In this test, the number of FOM runs during the hybrid model simulation is $N_f=469$ times, which takes about $P_f=22.9\%$ of the total number of time steps. 
The CPU time for integrating the hybrid model is $t_h=2.03$ seconds while the full order model simulation takes $t_f=4.25$ seconds. 
Overall, it is observed that the hybrid approach in this test is able to achieve accurate state and output approximations close to those of the FOM at a lower cost.  
 
Secondly, we investigate the numerical behavior of the hybrid model with respect to the output error tolerance and the width of time window. To this end, we fix the estimated error tolerance $\mathsf{tol}$ to be $10^{-3}$, $10^{-4}$ or $10^{-5}$ while varying the time window width $w$ from $10, 20$ to $40$. The results are listed in Table \ref{tab:Burg1dnu1}.  
It is seen that, for a fixed tolerance, varying $w$ does not have a significant influence on the output error and the number of full order simulations in the hybrid approach. 
Meanwhile, for a fixed $w$, decreasing the tolerance would reduce the output error, but resulting more FOM runs in the hybrid approach, which also cause more CPU time to complete the simulation. 
\begin{rem}\label{rem1}
Note that $\Delta^{k+1}$ estimates the output error at the $k$-th step. Because the actual error accumulates over time, the time average of output error $\mathcal{E}_o$ can be greater than the user-defined error tolerance at each time step, $\mathsf{tol}$. 
\end{rem}

\begin{table}[ht]
\centering 
\begin{tabular}{c | c c c | c c c | c c c} 
\hline\hline 
\multirow{2}{*}{w} & \multicolumn{3}{c|}{$\mathsf{tol}=10^{-3}$} & \multicolumn{3}{|c|}{$\mathsf{tol}=10^{-4}$} & \multicolumn{3}{|c}{$\mathsf{tol}=10^{-5}$} \\ 
\cline{2-10}
              {} & $\mathcal{E}_o$ & $P_f$ & $t_h$ & $\mathcal{E}_o$ & $P_f$ & $t_h$ & $\mathcal{E}_o$ & $P_f$ & $t_h$ \\ [0.5ex] 
\hline 
10 & $1.26\times 10^{-2}$ & 14.0\% & 1.30 & $8.87\times 10^{-4}$ & 22.8\% & 1.73 & $6.41\times 10^{-5}$ & 30.1\% & 2.20 \\ 
20 & $1.28\times 10^{-2}$ & 13.8\% & 1.40 & $8.79\times 10^{-4}$ & 22.9\% & 2.03 & $6.57\times 10^{-5}$ & 29.7\% & 2.49 \\
40 & $1.28\times 10^{-2}$ & 13.8\% & 1.77 & $9.04\times 10^{-4}$ & 22.5\% & 2.89 & $6.62\times 10^{-5}$ & 29.3\% & 2.93 \\ [1ex] 
\hline 
\end{tabular}
\caption{The hybrid model results at different choices of error tolerance $\mathsf{tol}$ and time window width $w$ for 1D Burgers' equation with $\nu=10^{-3}$.} 
\label{tab:Burg1dnu1} 
\end{table}

\begin{figure}[!htb]
 \centering
        \includegraphics[width=1\linewidth]{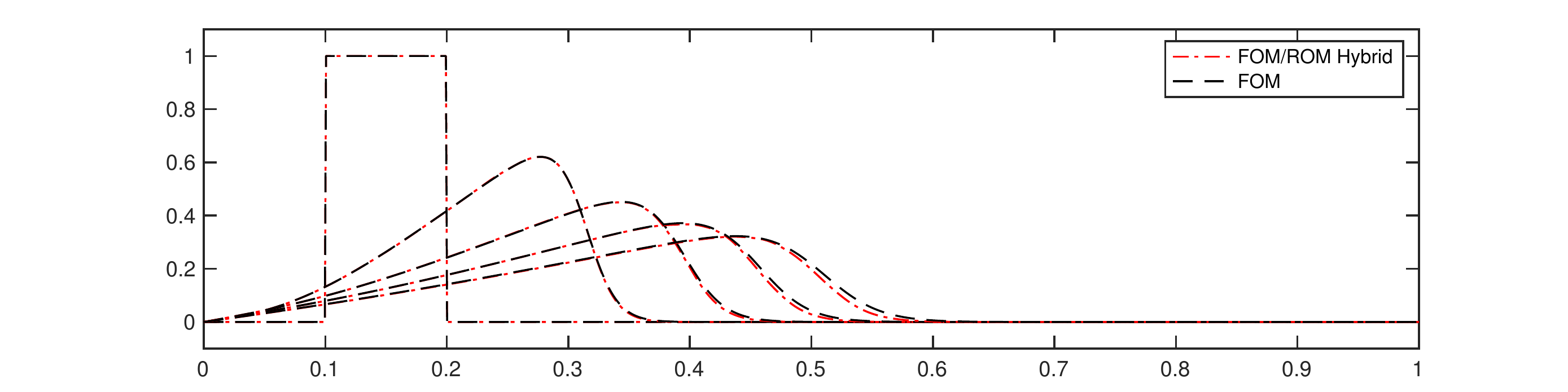}
\caption{Time snapshots $t= 0, 0.25, 0.5, 0.75, 1$ of Burgers' equation at $\nu=5\times10^{-3}$, $w=20$. \label{fig:snap2}}
\end{figure}

\begin{figure}[!htb]
 \centering
        \includegraphics[width=1\linewidth]{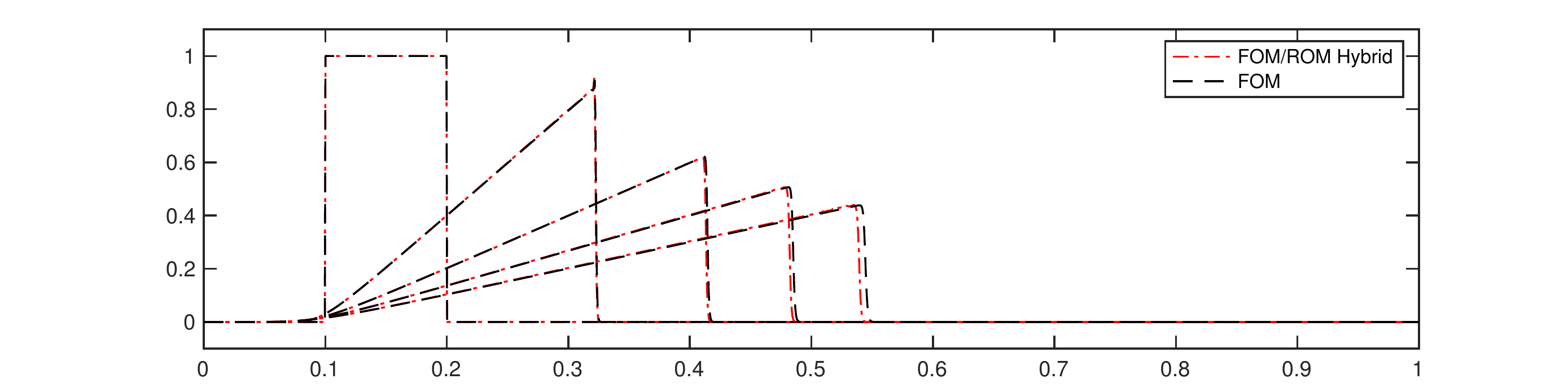}
\caption{Time snapshots $t= 0, 0.25, 0.5, 0.75, 1$ of Burgers' equation at $\nu=2\times10^{-4}$, $w=20$. \label{fig:snap3}}
\end{figure}

Thirdly, we fix the parameters $w=20$ and $\mathsf{tol}=10^{-4}$ in the hybrid approach but enlarge or shrink the viscosity coefficients of the Burger' equation by 5 times. 
When $\nu=5\times 10^{-3}$, we set $p=10$ in the discretization. 
The associated numerical results at selected time instances are plotted in Figure \ref{fig:snap2}.
When $\nu=2\times 10^{-4}$, we change $p=11$ because the problem is more convection dominated. 
Snapshots of numerical results at the same time instances are presented in Figure \ref{fig:snap3}.
The output error, percentage of FOM runs, and the CPU times for both cases are listed in Table \ref{tab:Burg1dnu2}, respectively. 
It is seen that the FOM runs are conducted at 10.4\% of the entire time steps in the first case and about 32.1\% in the second case.


\begin{table}[ht]
\centering 
\begin{tabular}{c c c c | c c c c } 
\hline\hline 
$\nu$ & $\mathcal{E}_o$ & $P_f$ & $t_h\,\, (t_f)$ & $\nu$ & $\mathcal{E}_o$ & $P_f$  & $t_h\,\,(t_f)$\\ [0.5ex] 
\hline 
$2\times 10^{-3}$ & $1.05\times 10^{-3}$ & 10.4 \% & 1.16 (4.07) & $2\times 10^{-4}$ & $8.35\times10^{-4}$ & 32.1\% & 27.2 (44.4) \\
\hline 
\end{tabular}
\caption{The hybrid model results for Burgers' equation when $w=20$ and $\mathsf{tol}=10^{-4}$.} 
\label{tab:Burg1dnu2} 
\end{table}

\paragraph{Test case 2.} We next consider the 1D parameterized Burgers's equation. The computational setting is the same as that in the first test case, except the diffusion coefficient $\nu\in [0.005, 1]$. To build a ROM valid in the parameter domain, we employ the hybrid model-based adaptive POD-Greedy-DEIM method presented in Algorithm \ref{alg:adp}. Although the nonlinearity in the equation is quadratic, which can be treated efficiently by tensor calculation, we use the DEIM approach here.  
We compare the performance of the adaptive basis selection algorithm when the hybrid model is used with that when the FOM is considered. 
 
For the offline training, $\Theta$ is composed of 20 logarithmically spaced parameter values in $[0.005, 1]$ and the tolerance is set to be $\mathsf{tol}=10^{-3}$. When the FOM is used as the high-fidelity solver for the snapshot generation, after adaptively selecting the POD and DEIM basis, it ends up with a final dimension $(\ell_{POD}, \ell_{DEIM})= (21,21)$. The whole process takes 7 iterations and involves 3 FOM runs. The CPU time for these 3 FOM runs takes 93.9 seconds in total. On the other hand, when the FOM/ROM hybrid model is used as the high-fidelity solver, after the greedy search, it ends up with a final dimensional $(\ell_{POD}, \ell_{DEIM}) = (25,26)$. The process takes 10 iterations and has 5  hybrid model runs. The CPU time for the hybrid model simulations is 30.6 seconds in total. The results are listed in Table \ref{tab:pBurg}. It is observed that using the hybrid model saves about 2/3 of offline time compared with the FOM used in the adaptive basis selection algorithm. 

\begin{table}[ht]
\centering 
\begin{tabular}{c | c c c c } 
\hline\hline 
 \multirow{ 2}{*}{FOM is used} & $(\ell_{POD}, \ell_{DEIM})$ & iterations & FOM & CPU of FOM \\ [0.5ex] 
\cline{2-5} 
  & $(21, 21)$ & 7 & 3 runs & 93.90 second \\
\hline 
  \multirow{ 2}{*}{Hybrid is used} & $(\ell_{POD}, \ell_{DEIM})$ & iterations & Hybrid & CPU of hybrid \\ [0.5ex] 
\cline{2-5} 
  & $(25, 26)$ & 10 & 5 runs & 30.63 second \\
\hline 
\end{tabular}
\caption{Adaptive basis selection when the FOM and the hybrid model are respectively used for snapshot generations.} 
\label{tab:pBurg} 
\end{table}
\begin{figure}[!htb]
 \centering
        \includegraphics[width=1\linewidth]{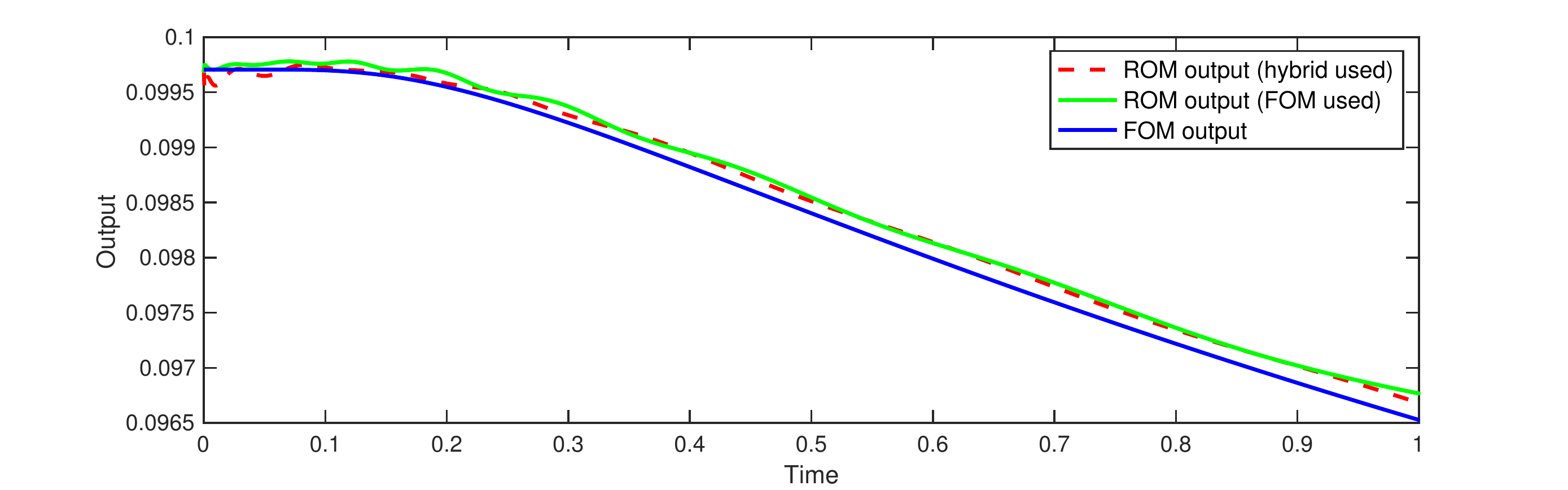}
\caption{Outputs generated from the FOM and the ROMs. Either FOM or hybrid model is used for adaptively selecting basis vectors. \label{fig:pBurg1}}
\end{figure}

To check the performance of the ROM, we consider the $\nu=0.005$ case. The time evolution of the output in the reduced-order approximations are plotted in Figure \ref{fig:pBurg1}, together with the output of the full-order simulation. It shows that the ROM generated by either the adaptive POD-Greedy-DEIM or the hybrid model-based adaptive POD-Greedy-DEIM algorithm is able to provide accurate outputs, which are close to the FOM outputs.   

\begin{figure}[!htb]
\centering
\includegraphics[width=0.9\linewidth, height=0.3\linewidth]{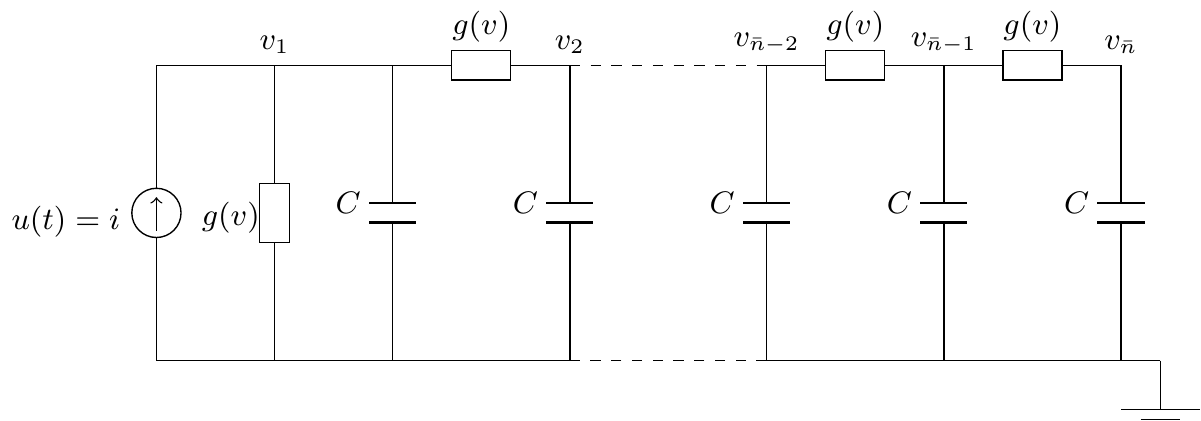}
\caption{RC circuit diagram.}
\label{figrc1}
\centering
\end{figure}

\paragraph{Test case 3.} Next we consider a nonlinear circuit example with RC structure shown in Figure~\ref{figrc1}. 
This is a widely used benchmark example for MOR of nonlinear systems~\cite{morGu11,morAsiABetal21}.
The diodes in the circuit constitute the nonlinear part of the model, which characterized by $g(v)=e^{40v}-1$ with $v$ the nodal voltage. The voltage $v_{1}(t)$ at node 1 is the output and the current $I$ is the input signal $u(t)$ of the system. After applying Kirchhoff's current law at each of the $\bar{n}$ nodes, and assuming that the  capacitance is normalized, i.e., $C=1$, we obtain the nonlinear circuit model
\begin{equation}
\dot{v}(t)=Av(t)+f(v(t))+bu(t),~~ y(t)=cv(t),
\end{equation}
where Taylor series expansion at $v=0$ is applied to separate the linear part $A v$ from $g(v)$ so that the proposed error indicator can be efficiently applied, i.e., $g(v)=Av(t)+f(v(t))$. The input, output matrices are given by $b=c^{T}\in \mathbb R^{\bar{n}}$, and $b$ is a unit basis vector with 1 in the first entry and all other entries are zeros. The input $u(t)$ is a step signal, taken as $u(t)=0$ if $t\leq 3$ and $u(t)=1$ if $3<t\leq 10$. 


\begin{figure}[!htb]
 \centering
        \includegraphics[width=1\linewidth,height=2in]{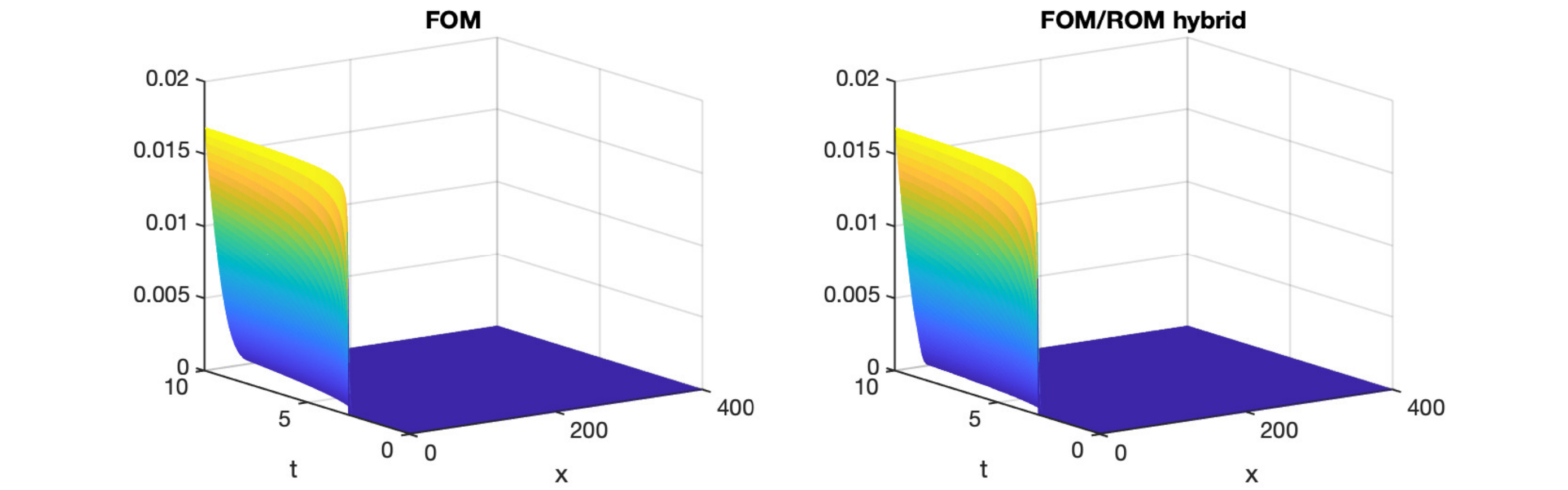}
\caption{Time evolution plots of nonlinear circuit with the step input signal: (left) FOM, (right) the hybrid model when $w=20$ and $\mathsf{tol}=10^{-4}$. \label{fig:evo1}}
\end{figure}
\begin{figure}[!htb]
 \centering
        \includegraphics[width=1\linewidth]{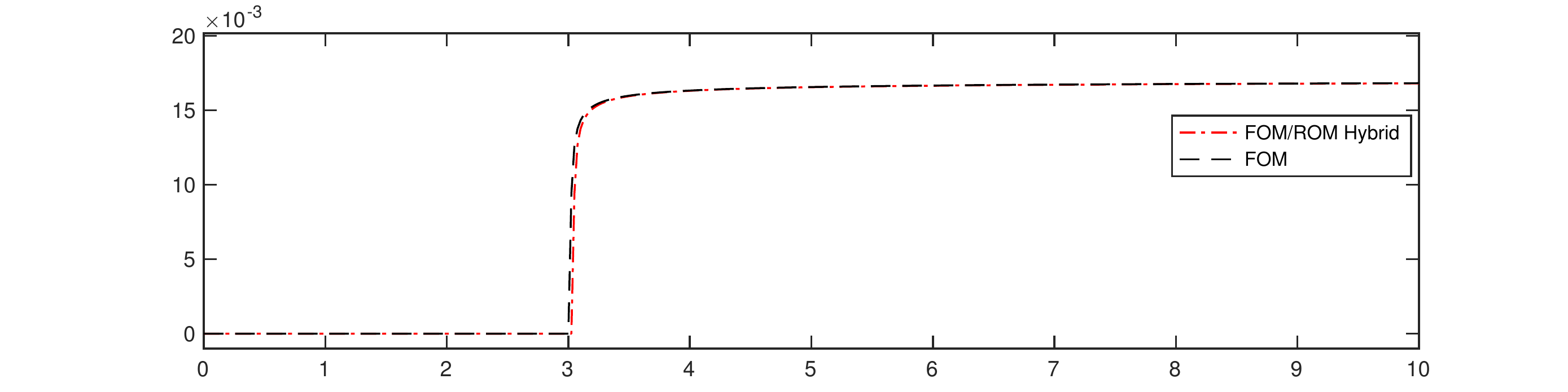}
\caption{Time evolution of the outputs from the FOM and the hybrid model. \label{fig:circuit_output}}
\end{figure}
\begin{figure}[!htb]
 \centering
        \includegraphics[width=1\linewidth]{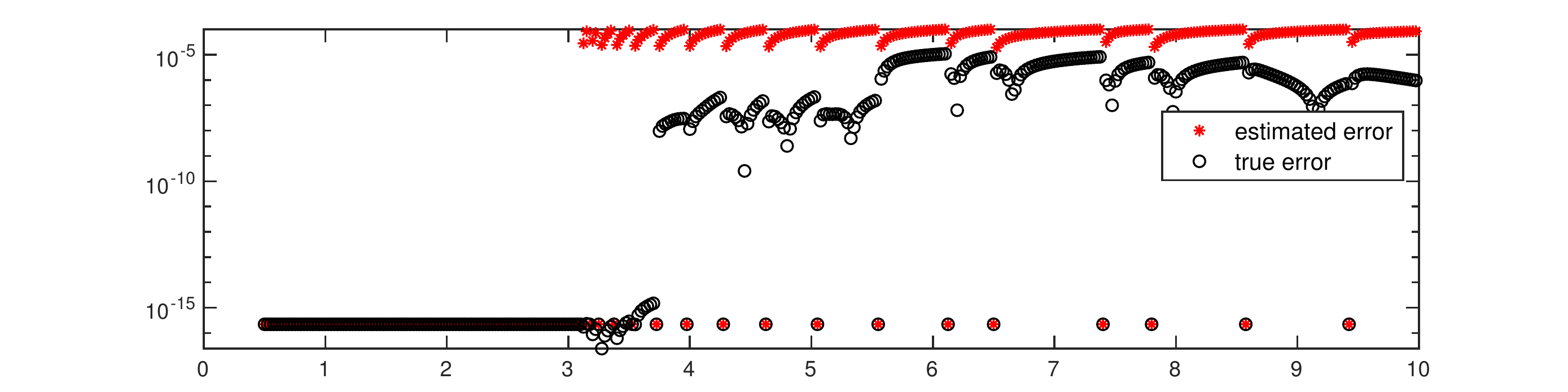}
\caption{Estimated and true errors of the circuit problem during the simulation.\label{fig:circuit_err}}
\end{figure}

We consider $\bar{n}= 401$ and the simulation time interval $t\in [0, 10]$ seconds in this test. For the time integration, we use semi-implicit method with $N_t= 400$ uniform time steps.   
In the hybrid approach, we use the POD ROM, set the time window width $w=20$ and the tolerance of output error at each step to be $10^{-4}$. 
The dual problem is solved once by a direct solver. 
The time evolutions of the FOM and the hybrid model voltage values are shown in Figure \ref{fig:evo1} and the output results are shown in Figure \ref{fig:circuit_output}, respectively. 
It is seen that, compared with the FOM, the hybrid approach obtains an accurate approximation. 
Time evolutions of the estimated error $\Delta^{k+1}$ and ``true" error are shown in Figure \ref{fig:circuit_err}, which indicates that $\Delta^{k+1}$ provides a reliable indicator of the output approximation error at each time step. 
Note that, because the system input is zero in the first 3 seconds, the estimated error is zero, which matches the real error in the time interval. To make these errors visible in Figure \ref{fig:circuit_err} which has a log-scale y-axis, we set them to be the floating-point relative accuracy ($\mathsf{eps}$ in Matlab notation). 

The total number of FOM runs during the hybrid model simulation is $N_f=40$ times, which includes $20$ initial FOM runs at the beginning of the simulation. The hybrid model takes $P_f=10\%$ of the total number of time steps. 
The CPU time for integrating the hybrid model is $t_h=4.22\times 10^{-2}$ seconds while the full order model simulation takes $t_f= 9.47\times 10^{-2}$ seconds. 

\section{Conclusions}
\label{sec:con}
Multiple high-fidelity large-scale simulations are generally needed in building ROMs for complex dynamical systems. 
These simulations provide snapshot data for which the reduced basis is extracted and the solution manifold is approximated.  
In this paper, a FOM/ROM hybrid approach is developed to accelerate such high-fidelity simulations. 
The development is based on an a posteriori error estimation for the output approximation. By controlling the estimated error, the method  switches between the FOM and the ROM generated on the fly in a dynamical manner. 
Numerical tests on the Burgers' equation and a nonlinear circuit problem illustrates that the hybrid approach is able to save the computational time while achieving a user-defined level of accuracy. Our future work includes the hybrid MOR for systems with changing input signals, the extension of the error indicator to state approximations and the design of more efficient algorithms for updating the ROM on the fly. 
\bibliographystyle{abbrv}
\bibliography{D_bibliography}
\end{document}